\input amstex
\input epsf
\documentstyle{amsppt}
\magnification 1200
\hsize=6truein
\vsize=8.5truein

\def\ls{\leqslant}
\def\gs{\geqslant}
\def\SS{\frak S}
\def\qq{{\bold q}}

\TagsOnRight

\topmatter
\title
Restricted permutations, continued fractions, and Chebyshev polynomials
\endtitle
\author Toufik Mansour$^*$ and Alek Vainshtein$^\dag$ \endauthor
\affil $^*$ Department of Mathematics\\
$^\dag$ Department of Mathematics and Department of Computer Science\\ 
University of Haifa, Haifa, Israel 31905\\ 
{\tt tmansur\@study.haifa.ac.il},
{\tt alek\@mathcs.haifa.ac.il}
\endaffil

\abstract
Let $f_n^r(k)$ be the number of 132-avoiding permutations on
$n$ letters that contain exactly $r$ occurrences of $12\dots k$,
and let $F_r(x;k)$ and $F(x,y;k)$ be the generating functions defined by
$F_r(x;k)=\sum_{n\gs0} f_n^r(k)x^n$ and $F(x,y;k)=\sum_{r\gs0}F_r(x;k)y^r$.
We find an explicit expression for $F(x,y;k)$ in the form of a continued 
fraction. This allows us to express $F_r(x;k)$ for $1\ls r\ls k$ via Chebyshev
polynomials of the second kind.
\medskip
\noindent {\smc 2000 Mathematics Subject Classification}: 
Primary 05A05, 05A15; Secondary 30B70, 42C05
\endabstract

\rightheadtext{Permutations, fractions, and Chebyshev polynomials}
\leftheadtext{Toufik Mansour and Alek Vainshtein}
\endtopmatter

\vfill\eject

\document
\heading 1. Introduction \endheading

Let $[p]=\{1,\dots,p\}$ denote a totally ordered alphabet on $p$ letters, and 
let $\alpha=(\alpha_1,\dots,\alpha_m)\in [p_1]^m$, 
$\beta=(\beta_1,\dots,\beta_m)\in [p_2]^m$. We say that $\alpha$ is
{\it order-isomorphic\/} to $\beta$ if for all $1\ls i<j\ls m$ one has
$\alpha_i<\alpha_j$ if and only if $\beta_i<\beta_j$. For two permutations
$\pi\in\SS_n$ and $\tau\in\SS_k$, an {\it occurrence\/} of $\tau$ in $\pi$
is a subsequence $1\ls i_1<i_2<\dots<i_k\ls n$ such that $(\pi_{i_1},
\dots,\pi_{i_k})$ is order-isomorphic to $\tau$; in such a context $\tau$ is
usually called the {\it pattern\/}. We say that $\pi$ {\it avoids\/} $\tau$,
or is $\tau$-{\it avoiding\/}, if there is no occurrence of $\tau$ in $\pi$. 
The set of all $\tau$-avoiding permutations
of all possible sizes including the empty permutation is denoted $\SS(\tau)$. 
 Pattern avoidance proved
to be a useful language in a variety of seemingly unrelated problems, from
stack sorting \cite{5} to singularities of Schubert varieties \cite{6}.
A complete study of pattern avoidance for the case $\tau\in\SS_3$ is
carried out in \cite{11}. For the case $\tau\in\SS_4$ see \cite{14, 
11, 12, 1}.

A natural generalization of pattern avoidance is the restricted pattern
inclusion, when a prescribed number of occurrences of $\tau$ in $\pi$ is 
required. Papers \cite{8} and \cite{3} contain simple expressions for
the number of permutations containing exactly one 123 and 132 patterns, 
respectively. The main result of \cite{B2} is that the generating function
for the number of permutations containing exactly $r$ 132 patterns is a
rational function in variables $x$ and $\sqrt{1-4x}$. This proves a
particular case of the general conjecture of Noonan and Zeilberger
\cite{9} which is that for any set $T$ of patterns, the sequence of numbers
enumerating permutations having a prescribed number of occurrences of 
patterns in $T$ is $P$-recursive. Recent paper \cite{10} presents the
generating function for the number of 132-avoiding permutations that contain 
a prescribed number of 123 patterns. The generating function is given in 
the form of a continued fraction. In the present  note we generalize the 
argument of \cite{10} to get the generating function for the number of 
132-avoiding permutations that contain a prescribed number of $12\dots k$ 
patterns for arbitrary $k\gs3$. The study of the obtained continued fraction
allows us to recover and to generalize the result of \cite{4} that relates
the  number of 132-avoiding permutations that contain no $12\dots k$ 
patterns to Chebyshev polynomials of the second kind.

The authors are grateful to C.~Krattenthaler, H.~Wilf, and
anonymous referee for useful comments
concerning Theorems~4.1 and~4.2. 

\heading 2. Continued fractions
\endheading

Let $f_n^r(k)$ stand for the number of 132-avoiding permutations on
$n$ letters that contain exactly $r$ occurrences of $12\dots k$. We denote
by $F(x,y;k)$ the generating function of the sequence $\{f_n^r(k\})$, that
is,
$$
F(x,y;k)=\sum_{n\gs0}\sum_{r\gs0}f_n^r(k)x^ny^r.
$$

Our first result is a natural generalization of the main theorem of
\cite{10}. 

\proclaim{Theorem 2.1} The generating function $F(x,y;k)$ for $k\gs1$ 
is given by the continued fraction
$$
F(x,y;k)=\frac1{1-\dfrac{xy^{d_1}}{1-\dfrac{xy^{d_2}}
{1-\dfrac{xy^{d_3}}{\dots}}}},
$$
where $d_i=\binom{i-1}{k-1}${\rm,} and $\binom ab$ is assumed $0$ whenever
$a<b$ or $b<0$.
\endproclaim

\demo{Proof} Following \cite{10} we define $\eta_j(\pi)$, $j\gs1$, as the
number of occurrences of $12\dots j$ in $\pi$. Define $\eta_0(\pi)=1$
for any $\pi$, which means that the empty pattern occurs exactly once in
each permutation. The {\it weight\/} of a permutation $\pi$ is a monomial
in $k$ independent variables $q_1,\dots,q_k$ defined by
$$
w_k(\pi)=\prod_{j=1}^k q_j^{\eta_j(\pi)}.
$$
The {\it total weight\/} is a polynomial
$$
W_k(q_1,\dots,q_k)=\sum_{\pi\in\SS(132)}w_k(\pi).
$$
The following proposition is implied immediately by the definitions.

\proclaim{Proposition 2.1} $F(x,y;k)=W_k(x,1,\dots,1,y)$ for $k\gs2${\rm,}
and $F(x,y;1)=W_1(xy)$.
\endproclaim

We now find a recurrence relation for the numbers $\eta_j(\pi)$.
Let $\pi\in\SS_n$, so that $\pi=(\pi',n,\pi'')$.

\proclaim{Proposition 2.2} For any $j\gs1$ and any nonempty $\pi\in\SS(132)$
$$
\eta_j(\pi)=\eta_j(\pi')+\eta_j(\pi'')+\eta_{j-1}(\pi').
$$ 
\endproclaim

\demo{Proof} Let $l=\pi^{-1}(n)$.
Since $\pi$ avoids 132, each number in $\pi'$ is greater than
any of the numbers in $\pi''$. Therefore, $\pi'$ is a 132-avoiding
permutation of the numbers $\{n-l+1,n-l+2,\dots,n-1\}$, while $\pi''$
is a 132-avoiding permutation of the numbers $\{1,2,\dots,n-l\}$. On the
other hand, if $\pi'$ is an arbitrary 132-avoiding
permutation of the numbers $\{n-l+1,n-l+2,\dots,n-1\}$ and $\pi''$ is an 
arbitrary 132-avoiding permutation of the numbers $\{1,2,\dots,n-l\}$, then
$\pi=(\pi',n,\pi'')$ is 132-avoiding. Finally, if $(i_1,\dots,i_j)$ is
an occurrence of $12\dots j$ in $\pi$ then either $i_j<l$, and so it is
also an occurrence  of $12\dots j$ in $\pi'$, or $i_1>l$, and so it is
also an occurrence  of $12\dots j$ in $\pi''$, or $i_j=l$, and so
$(i_1,\dots,i_{j-1})$ is an occurrence  of $12\dots j-1$ in $\pi'$.
The result follows.
\qed
\enddemo

Now we are able to find the recurrence relation for the total weight $W$.
Indeed, by Proposition~2.2, 
$$
\align
W_k(q_1,\dots,q_k)&=1+\sum_{\varnothing\ne\pi\in\SS(132)}\prod_{j=1}^k
q_j^{\eta_j(\pi')+\eta_j(\pi'')+\eta_{j-1}(\pi')}\\
&=1+\sum_{\pi'\in\SS(132)}\sum_{\pi''\in\SS(132)}\prod_{j=1}^k
q_j^{\eta_j(\pi'')}\cdot q_1\prod_{j=1}^{k-1}(q_jq_{j+1})^{\eta_j(\pi')}
\cdot q_k^{\eta_k(\pi')}\\
&=1+q_1W_k(q_1,\dots,q_k)W_k(q_1q_2, q_2q_3,\dots,q_{k-1}q_k,q_k).
\tag1\endalign
$$

For any $d\gs0$ and $1\ls m\ls k$ define
$$
\qq^{d,m}=\prod_{j=1}^k q_j^{\binom{d}{j-m}};
$$
recall that $\binom ab=0$ if $a<b$ or $b<0$. The following proposition is
implied immediately by the well-known properties of binomial coefficients.

\proclaim{Proposition 2.3} For any $d\gs0$ and $1\ls m\ls k$
$$
\qq^{d,m}\qq^{d,m+1}=\qq^{d+1,m}.
$$
\endproclaim

Observe now that $W_k(q_1,\dots,q_k)=W_k(\qq^{0,1},\dots,\qq^{0,k})$ and that
by (1) and Proposition~2.3
$$
W_k(\qq^{d,1},\dots,\qq^{d,k})=
1+\qq^{d,1}W_k(\qq^{d,1},\dots,\qq^{d,k})W_k(\qq^{d+1,1},\dots,\qq^{d+1,k}),
$$
therefore
$$
W_k(q_1,\dots,q_k)=\frac1{1-\dfrac{\qq^{0,1}}{1-\dfrac{\qq^{1,1}}
{1-\dfrac{\qq^{2,1}}{\dots}}}}.
$$
To obtain the continued fraction representation for $F(x,y;k)$ it is enough
to use Proposition~2.1 and to observe that
$$
\qq^{d,1}\bigg|_{q_1=x,q_2=\dots=q_{k-1}=1,q_k=y}
=xy^{\binom{d}{k-1}}. \qquad\qed
$$
\enddemo

\remark{Remark} For $k=1$ one recovers from Theorem~2.1 the well-known 
generating
function for the Catalan numbers, $(1-\sqrt{1-4z})/2z$. This result also
follows immediately from Proposition~2.1 and equation (1), which for $k=1$
is reduced to $W_1(q)=1+qW^2_1(q)$.
\endremark

\heading 3. Chebyshev polynomials
\endheading

Let us denote by $F_r(x;k)$ the generating function of the sequence 
$\{f_n^r(k)\}$ for a given $r$, that is,
$$
F_r(x;k)=\sum_{n\gs0} f_n^r(k)x^n.
$$
Recall that $F(x,y;k)=\sum_{r\gs0}F_r(x;k)y^r$. In this section we find
explicit expressions for $F_r(x;k)$ in the case $0\ls r\ls k$.

Consider a recurrence relation
$$
T_j=\frac1{1-xT_{j-1}},\quad j\gs1. \tag2
$$
The solution of (2) with the initial condition $T_0=0$ is denoted by
$R_j(x)$, and the solution of (2) with the initial condition
$$
T_0=G(x,y;k)=\frac y{1-\dfrac{xy^{\binom k1}}{1-\dfrac{xy^{\binom{k+1}2}}
{1-\dfrac{xy^{\binom{k+2}3}}{\dots}}}}
$$
is denoted by $S_j(x,y;k)$, or just $S_j$ when the value of $k$ is clear from
the context. Our interest in (2) is stipulated by the
following relation, which is an easy consequence of Theorem~2.1:
$$
F(x,y;k)=S_k(x,y;k). \tag3
$$

First of all, we find an explicit formula for the functions $R_j(x)$.
Let $U_j(\cos\theta)=\sin(j+1)\theta/\sin\theta$ be the Chebyshev polynomials
of the second kind.

\proclaim{Lemma 3.1} For any $j\gs1$
$$
R_j(x)=\frac{U_{j-1}\left(\frac1{2\sqrt{x}}\right)}
{\sqrt{x}U_{j}\left(\frac1{2\sqrt{x}}\right)}.
\tag4
$$
\endproclaim

\demo{Proof} Indeed, it follows immediately from (2) that $R_j(x)$ is the
$j$th approximant for the continued fraction
$$
\frac 1{1-\dfrac{x}{1-\dfrac{x}{1-\dfrac{x}{\dots}}}}.
$$
Hence, by \cite{7, Theorem~2, p.~194}, for any $j\gs1$ one has 
$R_j(x)=A_j(x)/A_{j+1}(x)$, where
$$
A_j(x)=\left(\frac{1+\sqrt{1-4x}}2\right)^j-
\left(\frac{1-\sqrt{1-4x}}2\right)^j.
$$
Using substitution $x\to 1/{4t^2}$ one gets 
$(2t)^jA_j(1/4t^2)=2\sqrt{t^2-1}U_{j-1}(t)$, which gives 
$A_j(x)=\sqrt{1/x-4}\,x^{j/2}U_{j-1}(1/2\sqrt{x})$, and the result follows.
\qed
\enddemo

Next, we find an explicit expression for $S_j$ in terms of $G$
and $R_j$.

\proclaim{Lemma 3.2} For any $j\gs1$ and any $k\gs1$
$$
S_j(x,y;k)=R_j(x)\frac{1-xR_{j-1}(x)G(x,y;k)}{1-xR_{j}(x)G(x,y;k)}. \tag5
$$
\endproclaim

\demo{Proof} Indeed, from (2) and $S_0=G$ we get 
$S_1=1/(1-xG)$. On the other hand, $R_0=0$, $R_1=1$, so (5)
holds for $j=1$. Now let $j>1$, then by induction
$$
S_j=\frac1{1-xS_{j-1}}=\frac1{1-xR_{j-1}}\cdot
\frac{1-xR_{j-1}G}{1-\dfrac{x(1-xR_{j-2})R_{j-1}G}
{1-xR_{j-1}}}.
$$
Relation (2) for $R_j$ and $R_{j-1}$ yields 
$(1-xR_{j-2})R_{j-1}=(1-xR_{j-1})R_j=1$, which together with the above
formula gives (5).
\qed
\enddemo

As a corollary from Lemma~3.2 and (3) we get the following expression for
the generating function $F(x,y;k)$.

\proclaim{Corollary}
$$
F(x,y;k)=R_k(x)+\big(R_k(x)-R_{k-1}(x)\big)
\sum_{m\gs1}\big(xR_k(x)G(x,y;k)\big)^m.
$$
\endproclaim

Now we are ready to express the generating functions $F_r(x;k)$, 
$0\ls r\ls k$, via Chebyshev polynomials.

\proclaim{Theorem 3.1} For any $k\gs1${\rm,} $F_r(x;k)$ is a rational function
given by
$$
\align
F_r(x;k)&=\frac{x^{\frac {r-1}2}U_{k-1}^{r-1}\left(\frac1{2\sqrt{x}}\right)}
{U_{k}^{r+1}\left(\frac1{2\sqrt{x}}\right)},\quad 1\ls r\ls k,\\
F_0(x;k)&=\frac{U_{k-1}\left(\frac1{2\sqrt{x}}\right)}
{\sqrt{x}U_{k}\left(\frac1{2\sqrt{x}}\right)},
\endalign
$$
where $U_j$ is the $j$th Chebyshev polynomial of the second kind.
\endproclaim

\demo{Proof} Observe that $G(x,y;k)=y+y^{k+1}P(x,y)$, 
so from Corollary we get
$$
F(x,y;k)=R_k(x)+\big(R_k(x)-R_{k-1}(x)\big)
\sum_{m=1}^k\big(xR_k(x)\big)^my^m+y^{k+1}P'(x,y),
$$
 where $P(x,y)$ and $P'(x,y)$ are formal power series.
To complete the proof, it suffices to use (4) together with the identity
$U_{n-1}^2(z)-U_n(z)U_{n-2}(z)=1$,
which follows easily from the trigonometric identity
$\sin^2n\theta-\sin^2\theta=\sin(n+1)\theta\sin(n-1)\theta$.
\qed
\enddemo

For the case $r=0$ this result was proved by a different method in \cite{4}.

\heading 4. Further results
\endheading

There are several ways to generalize the results of the previous sections. 
First, one can try to get exact formulas for $F_r(x;k)$ in the case $r>k$.
The method described in Section~3 allows, in principle, to obtain such 
formulas, though they become more and more complicated. For example, the
following theorem gives an explicit expression for $F_r(x;k)$ when
$r\ls k(k+3)/2$.

\proclaim{Theorem 4.1} For any $k\gs1$ and $1\ls r\ls k(k+3)/2$,
$F_r(x;k)$ is a rational function given by
$$
F_r(x;k)=\frac{x^{\frac {r-1}2}U_{k-1}^{r-1}\left(\frac1{2\sqrt{x}}\right)}
{U_{k}^{r+1}\left(\frac1{2\sqrt{x}}\right)}\sum_{j=0}^{\lfloor(r-1)/k\rfloor}
\binom{r-kj+j-1}{j}\left(\frac{U_k\left(\frac1{2\sqrt{x}}\right)}
{x^{\frac{k-2}{2k}}U_{k-1}\left(\frac1{2\sqrt{x}}\right)}\right)^{kj},
$$
where $U_j$ is the $j$th Chebyshev polynomial of the second kind. 
\endproclaim

\demo{Proof}
Indeed, the explicit expression for $G(x,y;k)$ gives
$$
G(x,y;k)=y(1+xy^k+\dots +x^sy^{ks})+y^tP(x,y),
$$
where $s=\lceil (k+1)/2\rceil$, $t=1+k(k+3)/2$, and $P(x,y)$ is a formal power
series. Hence, by Corollary,
$$\align
\frac{F(x,y;k)-R_k(x)}{R_k(x)-R_{k-1}(x)}&=
\sum_{m\gs1}\big(xR_k(x)\big)^my^m(1+xy^k+\dots +x^sy^{ks})^m+y^tP'(x,y)\\
&=\sum_{m\gs1}\big(xR_k(x)\big)^my^m
\sum_{j=0}^{ms}\binom{m+j-1}{j}x^jy^{kj}+y^tP'(x,y)\\
&=\sum_{r\gs1}y^r\big(xR_k(x)\big)^r\sum_{j=0}^{\lfloor(r-1)/k\rfloor}
\frac{\binom{r-kj+j-1}{j}x^j}{\big(xR_k(x)\big)^{kj}}
+y^tP''(x,y),
\endalign
$$
where $P'(x,y)$ and $P''(x,y)$ are formal power series. The rest
of the proof follows the proof of Theorem~3.1. 
\qed
\enddemo

Another possibility is to analyze the case of permutations containing
exactly one 132 pattern and $r$ $12\dots k$ patterns. Introducing
the modified total weight $\Omega_k(q_1,\dots,q_k)$ as the sum of the
weights $w_k(\pi)$ over all  permutations containing
exactly one 132 pattern, we get the following equation:
$$
\aligned
\Omega_k(q_1,\dots,q_k)&=
q_1W_k(q_1q_2,\dots, q_{k-1}q_k,q_k)\Omega_k(q_1,\dots,q_k)\\&+
q_1W_k(q_1,\dots,q_k)\Omega_k(q_1q_2,\dots, q_{k-1}q_k,q_k)\\&+
q_1^2q_2^2W_k(q_1q_2,\dots, q_{k-1}q_k,q_k)\big(W_k(q_1,\dots,q_k)-1\big);
\endaligned
$$
for the case $k=3$ see \cite{10}. By (1) and Proposition~2.3 this is
equivalent to 
$$
\aligned
\Omega_k(\qq^{d,1},\dots,\qq^{d,k})&=\qq^{d,1}\left(\qq^{d,2}\right)^2
\big(W_k(\qq^{d,1},\dots,\qq^{d,k})-1\big)^2\\
&+\qq^{d,1}W_k^2(\qq^{d,1},\dots,\qq^{d,k})
\Omega_k(\qq^{d+1,1},\dots,\qq^{d+1,k}).
\endaligned\tag6
$$
Let now $\varphi_n^r(k)$ be the number
of permutations on $n$ letters that contain exactly one 132 pattern and 
$r$ $12\dots k$ patterns, and $\Phi_r(x;k)$ be the generating function
of the sequence $\{\varphi_n^r(k)\}$ for a given $r$. In general, equation~(6)
allows us to find explicit expressions for $\Phi_r(x;k)$. However, they are
rather cumbersome, so we restrict ourselves to the case $r=0$.

\proclaim{Theorem 4.2} For any $k\gs3${\rm,} $\Phi_0(x;k)$ is a rational
function given by
$$\align
\Phi_0(x;k)&=\frac{x}{U_{k}^2\left(\frac1{2\sqrt{x}}\right)}
\sum_{j=1}^{k-2}U_{j}^2\left(\tfrac1{2\sqrt{x}}\right)\\
&=\frac 1{16\sin^2(k+1)t\cos^2t}
\left(2k-5+4\cos^2t-\frac{\sin(2k-1)t}{\sin t}\right),\endalign
$$
where $U_j$ is the $j$th Chebyshev polynomial of the second kind and
$\cos t=1/2\sqrt{x}$.
\endproclaim


\Refs
\ref 
\no 1
\by M.~Bona
\paper Permutations avoiding certain patterns: the case of length $4$ and some
generalizations
\jour Discr. Math.
\vol 175 \yr 1997 \pages 55--67
\endref

\ref 
\no 2
\by M.~Bona
\paper The number of permutations with exactly $r$ $132$-subsequences is
$P$-recursive in the size!
\jour Adv. Appl. Math.
\vol 18 \yr 1997 \pages 510--522
\endref

\ref 
\no 3
\by M.~Bona
\paper Permutations with one or two $132$-subsequences
\jour Discr. Math.
\vol 181 \yr 1998 \pages 267--274
\endref


\ref 
\no 4
\by T.~Chow and J.~West
\paper Forbidden subsequences and Chebyshev polynomials
\jour Discr. Math.
\vol 204 \yr 1999 \pages 119--128
\endref


\ref 
\no 5 
\by D.~Knuth
\book The Art of Computer Programming \vol 3
\publ Addison Wesley
\publaddr Reading, MA
\yr 1973
\endref


\ref 
\no 6
\by V.~Lakshmibai and B.~Sandhya
\paper Criterion for smoothness of Schubert varieties in 
$\operatorname{Sl}(n)/B$
\jour Proc. Indian Acad. Sci.
\vol 100 \issue 1 \yr 1990 \pages 45--52
\endref

\ref 
\no 7
\by L.~Lorentzen and H.~Waadeland
\book Continued fractions with applications
\publ North-Holland
\yr 1992
\endref


\ref 
\no 8
\by J.~Noonan
\paper The number of permutations containing exactly one increasing 
subsequence of length three
\jour Discr. Math.
\vol 152 \yr 1996 \pages 307--313
\endref

\ref 
\no 9
\by J.~Noonan and D.~Zeilberger
\paper The enumeration of permutations with a prescribed number of
``forbidden'' patterns
\jour Adv. Appl. Math.
\vol 17 \yr 1996 \pages 381--407
\endref

\ref 
\no 10
\by A.~Robertson, H.~Wilf, and D.~Zeilberger
\paper Permutation patterns and continuous fractions
\jour Elec. J. Comb.
\vol 6 \issue 1\yr 1999 \finalinfo R38 
\endref

\ref 
\no 11
\by R.~Simion and F.~Schmidt
\paper Restricted permutations
\jour Europ. J. Comb.
\vol 6\yr 1985 \pages 383-406
\endref

\ref 
\no 12
\by Z.~Stankova
\paper Classification of forbidden subsequences of length $4$
\jour Eur. J. Comb.
\vol 17 \yr 1996 \pages 501--517
\endref

\ref 
\no 13
\by Z.~Stankova
\paper Forbidden subsequences
\jour Discr. Math.
\vol 132 \yr 1994 \pages 291--316
\endref

\ref 
\no 14
\by J.~West
\paper Generating trees and the Catalan and Schr\"oder numbers
\jour Discr. Math.
\vol 146 \yr 1995 \pages 247--262
\endref
\endRefs
\enddocument